\def\vers{Apr.~8, 2009, v.1}
\magnification=1200
\hsize=6.5truein
\vsize=8.9truein
\font\bigfont=cmr10 at 14pt
\font\mfont=cmr9
\font\sfont=cmr8
\font\mbfont=cmbx9
\font\sifont=cmti8
\def\1{\hskip1pt}
\def\ssb{\raise.2ex\h{${\scriptscriptstyle\bullet}$}}
\def\msum{\h{$\sum$}}
\def\mopl{\h{$\bigoplus$}}
\def\a{\alpha}
\def\aa{\bar{a}}
\def\b{\beta}
\def\C{{\bf C}}
\def\CC{\bar{C}}
\def\D{\Delta}
\def\d{\partial}
\def\e{\eta}
\def\ee{\bar{\eta}}
\def\F{\hskip1pt\overline{\!F}\hskip1pt}
\def\H{{\bf H}}
\def\HH{{\bf H}^{(\lambda)}}
\def\h{\hbox}
\def\ha{\h{${1\over 2}$}}
\def\l{\langle}
\def\la{\lambda}
\def\lll{\bar{\lambda}}

\def\n{\nu}
\def\O{{\cal O}}
\def\Q{{\bf Q}}
\def\q{\quad}
\def\R{{\bf R}}
\def\r{\rangle}
\def\s{\sigma}
\def\t{\bar{t}}
\def\u{\widetilde{u}}
\def\uu{\bar{u}}
\def\V{{\cal V}}
\def\Vl{{\cal V}^{(\lambda)}}
\def\VV{\widehat{\cal V}}
\def\VVl{\widehat{\cal V}{}^{(\lambda)}}
\def\v{\widetilde{v}}
\def\vv{\bar{v}}
\def\w{\bar{w}}
\def\Z{{\bf Z}}
\def\z{\bar{z}}
\def\Gr{{\rm Gr}}
\def\Hom{{\rm Hom}}
\def\Ker{{\rm Ker}}
\def\Im{{\rm Im}}
\def\MHS{{\rm MHS}}
\def\NF{{\rm NF}}
\def\ad{{\rm ad}}

\def\too{\longrightarrow}
\def\into{\hookrightarrow}
\def\onto{\mathop{\rlap{$\to$}\hskip2pt\hbox{$\to$}}}
\def\simto{\buildrel\sim\over\too}
\def\({{\rm (}}
\def\){{\rm )}}
\h{}
\vskip 1cm
\centerline{\bigfont Cohomology classes of admissible normal
functions}

\bigskip
\centerline{Morihiko Saito}

\bigskip\medskip
{\narrower\noindent
{\mbfont Abstract.} {\mfont
We study the map associating the cohomology class of an admissible
normal function on the product of punctured disks, and give some
sufficient conditions for the surjectivity of the map. We also
construct some examples such that the map is not surjective.}
\par}

\bigskip\bigskip
\centerline{\bf Introduction}
\footnote{}{{\sifont Date\1}{\sfont:\ \vers}}

\bigskip\noindent
Let $\H$ be a polarizable variation of Hodge structure of weight
$-1$ on a product of punctured disks $S^*:=(\D^*)^n$.
Let $\NF(S^*,\H)_S^{\ad}$ be the group of admissible normal
functions with respect to $S=\D^n$ (see [12]).
Let $\NF(S^*,\H_{\Q})_S^{\ad}$ be its scalar extension by $\Z\to\Q$.
This is identified with the group of extension classes of
$\Q_{S^*}$ by $\H_{\Q}$ as admissible variations of $\Q$-mixed
Hodge structures with respect to $S$.
Let $j:S^*\to S$, $i_0:\{0\}\to S$ denote the inclusions.
We have a canonical morphism associating the cohomology class of an
admissible normal function
$$\NF(S^*,\H_{\Q})_S^{\ad}\to\Hom_{\MHS}(\Q,H^1i_0^*\R j_*\H_{\Q}),$$
where $\MHS$ denotes the category of graded-polarizable $\Q$-mixed
Hodge structures [6].
Let $j_{!*}\H_{\Q}$ denote the intermediate direct image [1].
This exists as a (shifted) pure Hodge module [11].
It is known (see e.g.\ [2] and also [10]) that
$H^1i_0^*j_{!*}\H_{\Q}$ is a subspace of $H^1i_0^*\R j_*\H_{\Q}$,
and the above morphism is naturally factored by
$$\NF(S^*,\H_{\Q})_S^{\ad}\to\Hom_{\MHS}(\Q,H^1i_0^*j_{!*}\H_{\Q}).
\leqno(0.1)$$
The target of (0.1) does not change by replacing $\H$ with
the nilpotent orbit associated to $\H$.
It is rather easy to show that (0.1) is surjective if $n=1$ (since
$H^1i_0^*j_{!*}\H_{\Q}=0$ in this case) or if $\H$ is a nilpotent
orbit or corresponds to a family of Abelian varieties.
More generally, we have

\medskip\noindent
{\bf Theorem~1.} {\it
Let $\V$ denote the underlying filtered $\O$-module of $\H$.
Assume $\Gr^{-2}_F\V=0$, or more generally, $F^{-1}\V$
is stable by the action of vector fields \(i.e.\ $F^{-1}\V$
is defined by a local system\). Then the morphism {\rm (0.1)} is
surjective.}

\medskip
For a more general statement, see Theorem~(1.5) below.
In this paper we also show that (0.1) is not necessarily surjective
for $n\ge 2$ in general.
More precisely, we have the following

\medskip\noindent
{\bf Theorem~2.} {\it Assume $\H$ is an nilpotent orbit of weight
$-1$ on $\D^*$ such that $\dim\Im\,N=1$ and $\Gr_F^1\V\ne 0$.
Then there is a \(non-horizontal\1\) one-parameter family of
polarizable variations of Hodge structures $\HH$ on $\D^*$
for $|\la|\ll 1$ \(shrinking $\D$ if necessary\1\) such that
$\H^{(0)}$ coincides with $\H$, the limit mixed Hodge structure
of $\HH$ is independent of $\la$, and moreover,
if $\HH$ denotes also its pull-back by the morphism
$S^*=(\D^*)^n\ni(t_1,\dots,t_n)\mapsto t_1\cdots t_n\in\D^*$
with $n\ge 2$, then for $0<|\la|\ll 1$, the morphism $(0.1)$
vanishes although its target does not, where
$H^1i_0^*j_{!*}\H_{\Q}\cong\mopl_{k=1}^{n-1}\Q$
as mixed Hodge structures.}

\medskip
Note that the situation in Theorem~2 is obtained by iterating
unnecessary blowing-ups along smooth centers contained in a smooth
divisor, and it is quite difficult to generalize it to a more
natural situation.
So Theorem~2 does not imply any obstructions to a strategy for
solving the Hodge conjecture which is recently studied by
M.~Green, P.~Griffiths, R.~Thomas, and others, see [2], [7], [15]
(and also Remark~(2.5) below).

In Section 1 we recall some basics of the cycle classes of
admissible normal functions, and prove Theorem~1.
In Section 2 we prove Theorem~2.

I would like to thank P.~Brosnan and G.~Pearlstein for useful
discussions and valuable comments.
This work is partially supported by Kakenhi 19540023.

\bigskip\bigskip
\centerline{\bf 1. Cohomology classes}

\bigskip\noindent
{\bf 1.1.~Calculation of cohomology classes.}
Let $S=\D^n$, $S^*=(\D^*)^n$, and $\H$ be a polarizable
variation of Hodge structure of weight $-1$ on $S^*$.
Let $H$ be the limit mixed Hodge structure of $\H$ with
$\Q$-coefficients, see [13].
We assume that the monodromies $T_i$ are unipotent.
Set $N_i={1\over 2\pi i}\log T_i$.
Let $j:S^*\to S$, $i_0:\{0\}\to S$ denote the inclusions.
Let $j_{!*}\H_{\Q}$ be the intermediate direct image [1].
Then $i_0^*j_{!*}\H_{\Q}$ is calculated by the complex
$$I^{\ssb}(H;N_1,\dots,N_n):=\bigl[0\to H\buildrel{\oplus_iN_i}
\over\too\mopl_i\,\Im\,N_i\to\mopl_{i\ne j}\,\Im\,N_iN_j\to
\cdots\bigr],\leqno(1.1.1)$$
where $H$ is put at the degree 0, see e.g. [5].
It is a subcomplex of the Koszul complex
$$K^{\ssb}(H;N_1,\dots,N_n):=\bigl[0\to H\buildrel{\oplus_iN_i}
\over\too\mopl_i\,H(-1)\to\mopl_{i\ne j}\,H(-2)\to\cdots\bigr],
\leqno(1.1.2)$$
calculating $i_0^*\R j_*\H_{\Q}$.
(Here $\Im\,N_i\subset H(-1)$, $\Im\,N_iN_j\subset H(-2)$,
see [6], 2.1.13 for Tate twist.)
This is obtained by iterating the restriction functors
$i_j^*=C(\psi_{t_j}\to\varphi_{t_j})$, where the $i_j$ denote the
inclusion of $\{t_j=0\}$ with $t_j$ the coordinates of the polydisk,
see [11].

We have a canonical morphism
$$\NF(S^*,\H)_S^{\ad}\to\Hom_{\MHS}(\Q,H^1i_0^*\R j_*\H_{\Q}),
\leqno(1.1.3)$$
associating the cohomology class of an admissible normal function
on $S^*$ with respect to $S$.
By definition [12], $\n\in\NF(S^*,\H)_S^{\ad}$ corresponds
(using [3]) to a short exact sequence of admissible variations of
Hodge structures ([8], [14])
$$0\to\H\to\H'\to\Z_{S^*}\to 0.\leqno(1.1.4)$$
Passing to the limit, we get a short exact sequence of mixed
$\Q$-Hodge structures endowed with the action of the logarithm of
the monodromies
$$0\to(H;N_1,\dots,N_n)\to(H';N'_1,\dots,N'_n)\to(\Q;0,\dots,0)
\to 0.$$
Let $\s:\Q\to H'$ be a splitting of the surjection $H'\to\Q$ as
$\Q$-vector spaces.
Then we have for any $i$
$$N'_i\s(1)\in\Im\,N_i,$$
restricting over curves transversal to each divisor
$\{t_i=0\}$, see e.g.\ [12], 2.5.4.
Hence
$$(N'_i\s(1))\in(\mopl_i\,\Im\,N_i)^0:=
\Ker(\mopl_i\,\Im\,N_i\to\mopl_{i\ne j}\,\Im\,N_iN_j),\leqno(1.1.5)$$
with
$$(\mopl_i\,\Im\,N_i)^0/\Im(\mopl_i\,N_i)=
H^1I^{\ssb}(H;N_1,\dots,N_n)=H^1i_0^*j_{!*}\H_{\Q}.$$

\medskip\noindent
{\bf 1.2.~Proposition.} {\it The cohomology class of an admissible
normal function $\n$ \(i.e.\ the image of $\n$ by $(1.1.3))$ is
given by}
$$(N'_i\s(1))\in H^1I^{\ssb}(H;N_1,\dots,N_n)=H^1i_0^*j_{!*}\H_{\Q}
\subset H^1i_0^*\R j_*\H_{\Q}.\leqno(1.2.1)$$

\medskip\noindent
{\it Proof.}
Applying $H^{\ssb}i_0^*\R j_*$ to (1.1.4), we get a long exact
sequence
$$\to H^0i_0^*\R j_*\Q_{S^*}\buildrel{\d}\over\to
H^1i_0^*\R j_*\H_{\Q}\to H^1i_0^*\R j_*\H'_{\Q}\to,\leqno(1.2.2)$$
and the cohomological class is given by the image of
$1\in\Q=H^0i_0^*\R j_*\Q_{S^*}$ by $\d$ (which is a morphism of
mixed Hodge structures).
Moreover, (1.2.2) is induced by the short exact sequence of
complexes of mixed Hodge structures
$$0\to K^{\ssb}(H;N_1,\dots,N_n)\to K^{\ssb}(H';N'_1,\dots,N'_n)\to
K^{\ssb}(\Q;0,\dots,0)\to 0,$$
where $K^{\ssb}(*)$ denotes the Koszul complex as in (1.1.2).
So the assertion follows.

\medskip\noindent
{\bf 1.3.~Mixed nilpotent orbits.}
We say that $((H,W');N_1,\dots,N_n)$ is a {\it mixed nilpotent
orbit} if $H$ is a mixed $\Q$-Hodge structure endowed with a
finite increasing filtration $W'$ and $N_i:H\to H(-1)$ are
nilpotent morphisms preserving $W'$ such that the following
two conditions are satisfied:

\medskip\noindent
(i) The relative monodromy filtration for $N_i$ with respect to
$W'$ exists for any $i$.

\smallskip\noindent
(ii) Each $(\Gr^{W'}_kH;N_1,\dots,N_n)$ is a pure nilpotent
orbit of weight $k$ for any $k$.

\medskip\noindent
Then the relative monodromy filtration for $\sum_{i\in I}N_i$
with respect to $W'$ exists for any subset $I$ of $\{1,\dots,n\}$,
see [8].
The category of mixed nilpotent orbits is an abelian category
such that any morphisms are strictly compatible with $F$ and $W'$,
see loc.~cit.
A mixed nilpotent orbit defines an admissible variation of mixed
Hodge structure on $S^*=(\D^*)^n$ with respect to $S=\D^n$ using
the coordinates $t_i$ of $\D^n$.
Here we use the correspondence between the multivalued horizontal
sections of $\V$ and the holomorphic s
ctions of the Deligne
extension $\VV$ of $\V$ annihilated by
$(t_k{\d\over\d t_k})^{\dim H}$ for any $k$.
It is defined by assigning to a horizontal section $v$
$$\v=\exp(-\msum_{k=1}^n(\log t_k)N_k)v.\leqno(1.3.1)$$

In the case $S^*=(\D^*)^n$, $S=\D^n$, and $\H$ is a nilpotent
orbit on $S^*$, we will denote by
$$\NF(S^*,\H_{\Q})_S^{\rm mno},
\leqno(1.3.2)$$
the subgroup of
$\NF(S^*,\H_{\Q})_S^{\ad}\,(=\NF(S^*,\H)_S^{\ad}\otimes_{\Z}\Q)$
consisting of admissible normal functions (tensored by $\Q$)
corresponding to extension classes in the category of mixed
nilpotent orbits.

\medskip
The following would be known to specialists.

\medskip\noindent
{\bf 1.4.~Proposition.} {\it Assume $\H$ is a nilpotent orbit.
Then $(0.1)$ is surjective.
More precisely, $(1.1.3)$ induces a surjective morphism}
$$\NF(S^*,\H_{\Q})_S^{\rm mno}\onto\Hom_{\MHS}(\Q,
H^1i_0^*j_{!*}\H_{\Q}).\leqno(1.4.1)$$

\medskip\noindent
{\it Proof.} It is enough to show the surjectivity of (1.4.1).
In the notation of (1.1.5), take
$$\a\in\Hom_{\MHS}(\Q,(\mopl_i\,\Im\,N_i)^0/\Im(\mopl_i\,N_i)).$$
This is identified with an element of
$(\mopl_i\,\Im\,N_i)^0/\Im(\mopl_i\,N_i)$ considering the image of
$1\in\Q$.
We have an exact sequence
$$0\to\Im(\mopl_i\,N_i)\to(\mopl_i\,\Im\,N_i)^0\to
(\mopl_i\,\Im\,N_i)^0/\Im(\mopl_i\,N_i)\to 0.$$
Take lifts $\a'_{\Q}$ and $\a'_F$ of $\a$ to
$(\mopl_i\,\Im\,N_i)^0_{\Q}$ and
$F^0(\mopl_i\,\Im\,N_i)^0_{\C}$ respectively.
There is $\b\in H_{\C}$ such that
$$(N_i(\b))=\a'_F-\a'_{\Q}\q\h{in}\,\,\,(\mopl_i\,\Im\,N_i)^0_{\C}.
\leqno(1.4.2)$$
We will construct an extension $H'$ of $\Q$ by $H$ such that the
image of the extension class by (1.4.1) corresponds to $\a$
as follows.

As a $\Q$-vector spaces we have
$$H'_{\Q}=H_{\Q}\oplus\Q.$$
The action of $N'_i$ on $H'_{\Q}$ is defined by
$$N'_i(a,b)=(N_ia+b(\a'_{\Q})_i,0)\q\h{for}\,\,\,a\in H_{\Q},\,
b\in\Q,\leqno(1.4.3)$$
where $(\a'_{\Q})_i\in(\Im\,N_i)_{\Q}$ is the $i$-th component of
$\a'_{\Q}$ in $(\mopl_i\,\Im\,N_i)^0_{\Q}$.
The weight filtration $W'$ is defined by $\Gr^{W'}_{-1}H'=H$
and $\Gr^{W'}_0H'=\Q$.
The Hodge filtration $F$ is defined by
$$F^pH'_{\C}=\cases{F^pH_{\C}&if $\,p>0$,\cr
F^pH_{\C}+\C(\b,1)&if $\,p\le 0$,\cr}\leqno(1.4.4)$$
where $F^pH_{\C}$ is identified with a subspace of $H'_{\C}$.

We have to show that $H'$ satisfies the conditions of mixed
nilpotent orbits.
By [8] it is enough to show that the relative monodromy
filtration exists for each $N'_i$, and the Griffiths
transversality $N'_iF^pH'_{\C}\subset F^{p-1}H'_{\C}$ is satisfied.
The first condition is trivially satisfied since
$(\a'_{\Q})_i\in(\Im\,N_i)_{\Q}$.
The second condition is reduced to
$$N'_i(\b,1)=N_i(\b)+(\a'_{\Q})_i\in F^{-1}H_{\C},$$
and follows from (1.4.2), i.e.\ $N_i(\b)+(\a'_{\Q})_i=(\a'_F)_i$.
(Note that the Hodge filtration $F$ on $\Im\,N_i\subset H_i(-1)$
is shifted by 1 so that $F^0(H(-1))_{\C}=F^{-1}H_{\C}$.)

By Proposition~(1.2) together with (1.4.3) for $(a,b)=(0,1)$,
the image by (1.4.1) of this extension class is given by
$\a'_{\Q}$.
So Proposition~(1.4) follows.

\medskip
The above argument can be extended as follows.
Let $\VV$ denote the Deligne extension of the underlying filtered
$\O$-module $\V$ of $\H$.
Then we have an isomorphism
$$\hbox{$H_{\C}\cong\bigcap_j\Ker\,(t_j{\d\over\d t_j})^{\dim H}
\subset\VV_0$,}$$
using the coordinates of $S$, and $H_{\C}$ is also identified with a
quotient $\V_0/m_0\V_0$ of $\V_0$ where $m_0$ denotes the maximal
ideal of $\O_{S,0}$.
Analyzing the proof of Proposition~(1.4), we get

\medskip\noindent
{\bf 1.5.~Theorem.} {\it
Assume the filtration induced by the inclusion $H_{\C}\into\VV_0$
coincides with the Hodge filtration of $H_{\C}$
which is by definition the quotient filtration by
$\VV_0\onto\VV_0/m_0=H_{\C}$,
i.e.\ $F^{-1}\VV_0$ is generated over $\O_{S,0}$ by its intersection
with $\bigcap_j\Ker\,(t_j{\d\over\d t_j})^{\dim H}\subset\VV_0$,
or equivalently, $F^{-1}\V$ is generated over $\O_{S^*}$ by the
sections $\v_1,\dots,\v_m$ of $\V$ corresponding to some horizontal
sections $v_1,\dots,v_m$ as in $(1.3.1)$.
Then $(0.1)$ is surjective}

\medskip\noindent
{\it Proof.} This follows from the same argument as in the proof of
Proposition~(1.4) if we replace (1.4.4) by
$$F^p\VV'=\cases{F^p\VV&if $\,p>0$,\cr
F^p\VV+\O_S(\b,1)&if $\,p\le 0$,\cr}$$
where $\V'$ is the underlying $\O$-module of the extension $\H'$,
and $\VV'$ is its Deligne extension.
Here $H'_{\C}$ is identified with a subspace of $\VV'$ using 
$\bigcap_j\Ker\,(t_j{\d\over\d t_j})^{\dim H'}$ as above.
The hypothesis implies that the image of $(\a'_F)_i\in H_{\C}$
in $\VV_0$ belongs to $F^{-1}\VV_0$ in the notation of the proof
of Proposition~(1.4), and hence the Griffiths transversality is
satisfied.
So Theorem~(1.5) follows.

\bigskip\bigskip
\centerline{\bf 2. Deformation of nilpotent orbits}

\bigskip\noindent
We will construct a family of variations of Hodge structures $\HH$
on $\D^*$ for $\la\in\C^*$ with $|\la|\ll 1$ by modifying the Hodge
filtration $F$ of $\H$ so that the limit mixed Hodge structure
does not change.
Here we may forget the rational structure and consider the real
structure instead, since we have the rational structure which
does not change by the deformation.

We will identify nilpotent orbit with the associated limit mixed
Hodge structure $H$ endowed with the action of $N$.
Let $\l *,*\r$ be a skew-symmetric form on $H_{\C}$ giving a
polarization of a nilpotent orbit as in [5].

\medskip\noindent
{\bf 2.1.~Proposition.} {\it With the assumption of Theorem~2,
there is a $4$-dimensional $\R$-nilpotent orbit $H_1$ generated by
$u_1,\dots,u_4\in H_{\C}$ satisfying the following conditions.

\smallskip\noindent
{\rm (i)} $u_1\in F^1W_{-1}H_{\C}$, $u_2\in W_0H_{\R}$,
$u_3\in W_{-2}H_{\C}$, $u_4\in F^{-2}W_{-1}H_{\C}$.

\smallskip\noindent
{\rm (ii)} $\uu_1=u_4$, $\uu_2=u_2$, $\uu_3=-u_3$, $Nu_2=u_3$,
$Nu_j=0\,\,(j\ne 2)$, where $N={1\over 2\pi i}\log T$.

\smallskip\noindent
{\rm (iii)} $[u_1],[u_4]\in\Gr^W_{-1}H_{\C}$ have type
$(1,-2)$ and $(-2,1)$ respectively.

\smallskip\noindent
{\rm (iv)} $u_2+au_3\in F^0W_0H_{\C}$ for some $a\in\C$.

\smallskip\noindent
{\rm (v)} $\l u_i,u_j\r=0$ unless $\{i,j\}=\{1,4\}$ or $\{2,3\}$.}

\medskip\noindent
{\it Proof.}
Since $\dim\Im\,N=1$, we have $N^2=0$ and the weight filtration $W$
on $H_{\R}$ is given by
$$W_{-3}H_{\R}=0,\q W_{-2}H_{\R}=\Im\,N,\q W_{-1}H_{\R}=\Ker\,N,\q
W_0H_{\R}=H_{\R}.$$
Moreover, we have an isomorphism compatible with $F$
$$N:\Gr^W_0H_{\R}\,(=\R)\simto(\Gr^W_{-2}H_{\R})(-1)\,(=(\R(1))(-1)).
\leqno(2.1.1)$$

By the second assumption of Theorem~2, there are
$$u_1\in F^1W_{-1}H_{\C},\q u_4\in F^{-2}W_{-1}H_{\C},$$
such that $\uu_1=u_4$ and condition (iii) is satisfied.
By the first assumption of Theorem~2, there is
$$u_3\in F^{-1}W_{-2}H_{\C},$$
which is purely imaginary, i.e. $\uu_3=-u_3$.
By (2.1.1) there is
$$u_2\in W_0H_{\R}\q\h{such that}\,\,\,Nu_2=u_3.$$
Note that
$$\Gr^W_0H_{\R}=\R[u_2],\q
\Gr^W_{-2}H_{\R}=\R[iu_3].\leqno(2.1.2)$$

To show condition (iv), note that $\Gr^W_{-1}H_{\C}$ is the direct
sum of $F^0$ and $\F^0$.
Since $i[v-\vv]\in\Gr^W_{-1}H_{\C}$ is real, it is written as
$[w]+[\w]$ for some $w\in F^0W_{-1}H_{\C}$.
Thus we may assume $v-\vv\in W_{-2}H_{\C}$ replacing $v$ with $v+iw$.
Then $v-au_3$ is real for some $a\in\C$, and we may replace
$u_2$ with $v-au_3$.

As for condition (v), it is satisfied modifying $u_1,u_4$ by a
multiple of $u_3$ if necessary, since $\l u_2,u_3\r\ne 0$ by
(2.1.2) and $\l u_j,u_3\r=0$ for $j=1,4$
(because $\l W_{-1},W_{-2}\r=0$).
This completes the proof of Proposition~(2.1).

\medskip\noindent
{\bf 2.2.~Construction.}
With the notation of Proposition~(2.1), let $H_1$ be the mixed
$\R$-Hodge structure spanned by $u_1,\dots,u_4$ in $H_{\C}$, and
$H_2$ be its orthogonal complement in $H_{\C}$.
We have the decomposition as mixed $\R$-Hodge structures endowed
with a pairing and the action of $N$
$$H_{\R}=H_1\oplus H_2.$$
For the proof of Theorem~2, we may then assume
$$H_{\R}=H_1,\q H_2=0,$$
so that $u_1,\dots,u_4$ is a $\C$-basis of $H_{\C}$.
Note that $H_{\C}$ is identified with the vector space of horizontal
sections of $\V$ and also with that of holomorphic sections of
the Deligne extension $\VV$ annihilated by
$(t_k{\d\over\d t_k})^{\dim H}$ for any $k$ (using (1.3.1)).

Let $z=\log t$ with $t$ the coordinate of $\D$.
The $u_j\,\,(j=1,\dots,4)$ induce a basis of the Deligne extension
$\VV$ as in (1.3.1), i.e.
$$\u_2=u_2-zu_3,\q\u_j=u_j\,\,(j\ne 2).$$
Setting $\xi=t{d\over dt}$, we have
$$\xi\u_2=-\u_3,\q \xi\u_j=0\,\,(j\ne 2).$$

Let $a,C\in\C$.
For $\la\in\C$ with $|\la|$ sufficiently small, define
$$\eqalign{&w_1=\u_1+C\la t(\u_2+a\u_3+\ha\la t\u_4),\cr
&w_2=\u_2+(a-1)\u_3+\la t\u_4,\q
w_3=-\u_3+\la t\u_4,\q w_4=\u_4.\cr}
\leqno(2.2.1)$$
Let $\VVl$ denote the Deligne extension of the underlying
$\O$-module $\Vl$ of $\HH$.
Define the Hodge filtration $F$ on $\VVl$ by
$$F^p=\msum_{i=1}^{2-p}\,\O_Sw_i.$$
For $\la=0$, the Hodge filtration on $\H^{(0)}$ coincides with
that of $\H$ choosing $a$ appropriately.
The Griffiths transversality holds for $\xi=t{d\over dt}$ since
$$\xi w_1=C\la tw_2,\q \xi w_2=w_3,\q \xi w_3=\la tw_4,\q
\xi w_4=0.\leqno(2.2.2)$$
However, this does not hold for $\la{\d\over \d\la}$,
and $\{\HH\}$ is a non-horizontal family.
Note that $Nu_2=u_3$, $Nu_j=0\,\,(j\ne 0)$, and
$$\eqalign{&w_1=u_1+C\la t(u_2-(z-a)u_3+\ha\la tu_4),\cr
&w_2=u_2-(z-a+1)u_3+\la tu_4,\q
w_3=-u_3+\la tu_4,\q w_4=u_4.\cr
&\w_1=u_4+\CC\lll\t(u_2+(\z-\aa)u_3+\ha\lll\t u_1),\cr
&\w_2=u_2+(\z-\aa+1)u_3+\lll\t u_1,\q
\w_3=u_3+\lll\t u_1,\q \w_4=u_1.\cr}\leqno(2.2.3)$$

Assume $C\in\C$ satisfies the condition
$$\l u_1,u_4\r=C\l u_2,u_3\r.$$
Then we have the orthogonal relation
$$\l F^p\Vl,F^{-p}\Vl\r=0\q\h{for any}\,\,\,p.
\leqno(2.2.4)$$
Indeed, the pairing is skew-symmetric and the above condition
on $C$ implies
$$\l w_1,w_2\r=\l w_1,w_3\r=0.$$

So it remains to show that
$F^p\cap\F{}^{-1-p}$ is 1-dimensional at each point of $\D^*$, and
for a generator $\e_{2-p}$ of $F^p\cap\F{}^{-1-p}$, we have the
positivity (as in [6], 2.1.15):
$$(2\pi i)^{-1}\l\e_{2-p},i^{-2p-1}\ee_{2-p}\r>0\,\,\,\,
(p=-2,-1,0,1).$$
It is enough to show these for $p=1$ and $0$, using the complex
conjugation.
For $p=1$ the assertion is easy (using $|\la|\ll 1$) since $F^1$
is generated by $w_1$.
For $p=0$, we first see that $F^0\cap\F{}^{-1}$ is at most
1-dimensional.
Indeed, otherwise it must coincide with $F^0$ and hence $\F^{-1}$
must contain $w_1$.
However, $\w_1,\w_2,\w_3,w_1$, or equivalently, $w_1,w_2,w_3,\w_1$
are linearly independent at each point of $\D^*$ if $|\la|\ll 1$.
The last assertion follows from (2.2.3) by calculating the
following determinant:
$$\pmatrix{1 & C\la t & -C\la t(z-a) & \ha C\la^2t^2\cr
0 & 1 & -z+a-1 & \la t\cr
0 & 0 & -1 & \la t\cr
\ha\CC\lll^2\t^2 & \CC\lll\t & \CC\lll\t(\z-\aa) & 1\cr}$$
Here we assume $\Im\,z\in[0,2\pi)$ so that $|tz|=|t\log t|$ is
bounded on $\D^*$.
Thus we get $\dim F^0\cap\F{}^{-1}\le 1$.
From (2.2.3) we also deduce
$$\eqalign{&w_1-C\la tw_2=u_1+C\la t(u_3-\ha\la tu_4)=
C\la t\w_3+(1-C\la\lll t\t)u_1-\ha C\la^2t^2u_4,\cr
&\w_1-\CC\lll\t\w_2=u_4+\CC\lll\t(-u_3-\ha\lll\t u_1)=
-\CC\lll\t\w_3+\ha\CC\lll^2\t^2u_1+u_4.\cr}$$
This implies that $u_1,u_4$ are linear combinations of
$w_1,w_2,\w_1,\w_2,\w_3$ with coefficients in
$\Q[\la,\lll,t,\t,P^{-1}]$
where $P=1-C\la\lll t\t+{1\over 4}C\CC\la^2\lll^2t^2\t{}^2$.
(Here $\la$ is viewed as a variable, but $a,c$ are constant.)
Substitute these to the following equality which also follows from
(2.2.3):
$$\eqalign{\w_2-w_2&=(z+\z-a-\aa+2)u_3+\lll\t u_1-\la tu_4\cr
&=(z+\z-a-\aa+2)\w_3-(z+\z-a-\aa+1)\lll\t u_1-\la tu_4.\cr}$$
Then we get
$$\la f_1w_1+(\la f_2+1)w_2=g_1\w_1+g_2\w_2+g_3\w_3,$$
where $f_1,f_2\in\Q[\la,\lll,t,\t,z\t,\z\t,P^{-1}]$ and
$g_1,g_2,g_3\in\Q[\la,\lll,t,\t,z,\z,P^{-1}]$.
Let $\e_2$ denote the left-hand side of the above equality so that
$$\e_2\in F^0\cap\F{}^{-1},\q \e_2-w_2=\la(f_1w_1+f_2w_2).$$
Since the pairing gives a polarization for $\la=0$ and
$|tz|=|t\log t|$ is bounded on $\D^*$
(assuming $\Im\,z\in[0,2\pi))$, the assertion follows.

\medskip\noindent
{\bf 2.3.~Theorem.} {\it
Let $\HH$ denote also the pull-back of the variation of Hodge
structure $\HH$ in $(2.2)$ by the morphism
$S^*:=(\D^*)^n\ni(t_1,\dots,t_n)\mapsto t_1\cdots t_n\in\D^*$.
Set $S=\D^n$.
Assume $\la\ne 0$ and $|\la|\ll 1$.
Then the morphism $(0.1)$ for $\HH$ vanishes although its
target does not, where $H^1i_0^*j_{!*}\HH_{\R}\cong
\mopl_{k=1}^{n-1}\R$ as mixed Hodge structures.}

\medskip\noindent
{\it Proof.}
By Proposition~(1.4) it is enough to show the vanishing of (0.1).
Since the assertion for the rational coefficients follows from
that for the real coefficients, we may assume that $H=H_1$ as in
(2.2).
We denote the pull-backs of $u_j,w_j\,\,(j=1,\dots,4)$
in (2.2) also by the same symbols.
Then $N_ku_j=0\,\,(j\ne 2)$ and $N_ku_2=u_3$ for $k=1,\dots,n$,
where $N_k={1\over 2\pi i}\log T_k$.
Since $N_k=N_{k'}$ for any $k,k'$, the $N_k$ will be denoted by $N$
(which can be viewed as the pull-back of $N$ on $\D^*$).
Let $z_k=\log t_k$.
The basis of the Deligne extension is given as in (1.3.1), i.e.
$$\u_2=u_2-\msum_{k=1}^n\,z_ku_3,\q\u_j=u_j\,\,(j\ne 2).
\leqno(2.3.1)$$
Setting $\xi_k=t_k{\d\over\d t_k}$, we have
$\xi_k\u_2=-\u_3$, $\xi_k\u_j=0\,\,(j\ne 2)$, and
$$\xi_kw_1=C\la t_1\cdots t_nw_2,\q\xi_kw_2=w_3,
\q\xi_kw_3=\la t_1\cdots t_nw_4,\q\xi_kw_4=0.\leqno(2.3.2)$$
Note that
$$\eqalign{&w_1=u_1+C\la t_1\cdots t_n(\u_2+au_3+\ha\la tu_4),\cr
&w_2=\u_2+(a-1)u_3+\la t_1\cdots t_nu_4,\q
w_3=-u_3+\la t_1\cdots t_nu_4,\q w_4=u_4.}$$

By the calculation in (1.1), the morphism (0.1) with real
coefficients is expressed as
$$\NF(S^*,\HH)^{\ad}_S\to H^1i_0^*j_{!*}\HH_{\R}=
(\mopl_{k=1}^n\,\Im\,N)/\Im\,N,\leqno(2.3.3)$$
since $\Im\,N\cong\R$ as mixed Hodge structures.
(Indeed, $H$ and $\Im\,N\subset H(-1)$ can be identified with those
associated to the variation on $\D^*$.)
Here $\Im\,N\into\mopl_{k=1}^n\,\Im\,N$ is the diagonal.
Thus the target of (2.3.3) is nonzero.
We have to show that (2.3.3) vanishes.

Assume there is an admissible normal function $\n$ such
that its image by (2.3.3) does not vanish.
Consider the corresponding short exact sequence of admissible
variations of $\R$-mixed Hodge structures
$$0\to\HH\to\HH{}'\to\R_{S^*}\to 0.$$
There is a basis $u'_j\,(j=0,\dots,4)$ of complex-valued horizontal
multivalued sections of $\Vl{}'$ such that
$u'_j=u_j\,\,(j\ne 0)$, $u'_0$ is defined over $\R$, and the image
of $u'_0$ in $\C$ is 1.
By (1.1) we have
$$N_ku'_0=c_ku_3\q\h{with}\,\,\,c_k\in\R\,\,(k=1,\dots,n),$$
and the image of the normal function $\n$ by (2.3.3) is
given by $(c_1,\dots,c_n)$, see (1.2.1).
So the above hypothesis on the image of $\n$ by (2.3.3) is
equivalent to
$$c_k\ne c_{k'}\q\h{for some}\,\,k,k'.\leqno(2.3.4)$$
Let $\u'_j$ be associated with $u'_j$ as in (2.3.1).
Then $\u'_j=\u_j\,\,(j\ne 0)$ and
$$\u'_0=u'_0-\msum_{k=1}^n\,z_kc_ku_3.$$
So we get for any $k$
$$\xi_k\u'_0=-c_k\u_3.\leqno(2.3.5)$$
Set $w'_j=w_j$ for $j=1,...,4$.
There are $h_1,\dots,h_4\in\O_S$ such that
$$w'_0:=\u'_0+\msum_{j=1}^4\,h_jw_j\in F^0\VVl,$$
(shrinking $S$ if necessary).
Then the Hodge filtration $F^p$ on the Deligne extension $\VVl$
is generated over $\O_S$ by $w'_i$ with $1\le i\le 2-p$
if $p>0$, and by $w'_0$ and $w'_i$ with $1\le i\le 2-p$ if $p\le 0$.
By the Griffiths transversality we have for any $k$
$$\xi_kw'_0\in F^{-1}\VVl=\msum_{j=0}^3\,\O_Sw'_j.$$
Using (2.3.2), (2.3.5) and the relation
$\u_3=\la t_1\cdots t_nw_4$ mod $\O_Sw_3$, we get for any $k$
$$\eqalign{(&\la t_1\cdots t_n(h_3-c_k)+\xi_kh_4)w_4\in F^{-1}\VVl,
\cr\h{i.e.}\q&\la t_1\cdots t_n(h_3-c_k)+\xi_kh_4=0.\cr}$$
This contradicts the hypothesis $c_k\ne c_{k'}$ in (2.3.4),
looking at the coefficient of $t_1\cdots t_n$ in the power series
expansion of $h_4$.
So the assertion follows.

\medskip\noindent
{\bf 2.4.~Remarks.}
(i) According to G.~Pearlstein, it is possible to describe the
examples as above using the theory of period maps as in
 [9], Th.~(6.16) or [4] Th.~(2.7).

\medskip
(ii) The arguments in (2.1--3) cannot be extended to arbitrary cases.
For example, we have the surjectivity of (0.1) if the hypothesis of
Theorem~1 or Theorem~(1.5) is satisfied.
We have the vanishing of the target of (0.1) if
$$F^{-1}\Gr^W_{-2}H_{\Q}\cap\Ker\,\Gr^WN=0.$$
Indeed, thia is shown by taking $\Gr^W$ of the differential
$$d:\mopl_i\,\Im\,N_i\to\mopl_{i\ne j}\,\Im\,N_iN_j$$
of the complex $I^{\ssb}(H;N_1,\dots,N_n)$ in (1.1.1).
Here $N_i=N$ for any $i$.

\medskip
(iii) In order to extend the arguments in (2.1--3),
it may be convenient to assume the following conditions
(which are stronger than the conditions coming from Remark~(ii)
above)
$$\Gr_F^{-2}(\Ker\,N)_{\C}\ne 0,\q
F^{-1}W_{-2}(\Ker\,N)_{\Q}\ne 0,\leqno(2.4.1)$$
where $\Ker\,N\subset H$ denotes a mixed Hodge structure.
(Indeed, if we do not assume the condition related to $\Ker\,N$,
then it does not seem easy to construct the Hodge filtration $F$
satisfying the Griffiths transversality (2.2.2) and the
orthogonality (2.2.4).)
It may be more convenient to replace the first condition of (2.4.1)
by a stronger one
$$\Gr_F^{-2}\Gr^W_{-1}(\Ker\,N)_{\C}\ne 0.$$

\medskip\noindent
{\bf 2.5.~Remark.}
As a consequence of recent work of R.~Thomas (see [14]), M.~Green
and P.~Griffiths (see [7]) and P.~Brosnan, G.~Pearlstein et al.\
(see [2]), it is
known that the Hodge conjecture is reduced to the non-vanishing
of the image by (0.1) of the normal function associated to a Hodge
class.
Theorem~2 does not imply any obstructions to this strategy since
the hypothesis of Theorem~2 is too strong.
It is very difficult to generalize the construction in (2.2) to
the situation appearing in [7] since the hypothesis
$\dim\Gr^W_{-2}H_{\Q}=\dim\Gr^W_0H_{\Q}=1$ is quite essential in
the argument of (2.2).
(This condition is essentially equivalent to that the divisor with
normal crossings is obtained by iterating unnecessary blowing-ups
along smooth centers contained in a smooth divisor.)

Assume, for simplicity, $n=3$ and
$$\Gr^W_{-2}H_{\Q}=\Q(1)\oplus\Q(1),\q\Gr^W_0H_{\Q}=\Q\oplus\Q,$$
and moreover $N_1,N_2,N_3$ are respectively expressed by
$$\pmatrix{1&0\cr0&0\cr},\pmatrix{0&0\cr0&1\cr},
\pmatrix{1&1\cr1&1\cr}.$$
Then it is quite difficult to construct a deformation of
the nilpotent orbit such that the image of (0.1) vanishes even
in this simple case.

\bigskip\bigskip
\centerline{{\bf References}}

\medskip
{\mfont
\item{[1]}
A.~Beilinson, J.~Bernstein and P.~Deligne, Faisceaux pervers,
Ast\'erisque, vol. 100, Soc. Math. France, Paris, 1982.

\item{[2]}
P.~Brosnan, H.~Fang, Z.~Nie and G.~J.~Pearlstein,
Singularities of admissible normal functions
(arXiv:0711.0964).

\item{[3]}
J.~Carlson, Extensions of mixed Hodge structures, in Journ\'ees
de G\'eom\'etrie Alg\'ebrique d'Angers 1979, Sijthoff-Noordhoff
Alphen a/d Rijn, 1980, pp. 107--128.

\item{[4]}
E.~Cattani and J.~Fernandez, Asymptotic Hodge theory and quantum
products, Contemp.\ Math., 276, Amer. Math. Soc., Providence, RI,
2001, pp.\ 115--136.

\item{[5]}
E.~Cattani, A.~Kaplan, W.~Schmid, $L^2$ and intersection
cohomologies for a polarizable variation of Hodge structure,
Inv.\ Math.\ 87 (1987), 212--252.

\item{[6]}
P.~Deligne, Th\'eorie de Hodge II, Publ. Math. IHES,
40 (1971), 5--57.

\item{[7]}
M.~Green and P.~Griffiths, Algebraic cycles and singularities
of normal functions, in Algebraic Cycles and Motives, Vol.~1
(J.~Nagel and C.~Peters, eds) London Math.\ Soc.\ Lect.\ Note
Series 343, Cambridge University Press 2007, 206--263.

\item{[8]}
M.~Kashiwara, A study of variation of mixed Hodge structure,
Publ.\ RIMS,\ Kyoto Univ. 22 (1986), 991--1024.

\item{[9]}
G.~J.~Pearlstein, Variations of mixed Hodge structure, Higgs
fields, and quantum cohomology, Manuscripta Math.\ 102 (2000),
269--310.

\item{[10]}
C.~Peters and M.~Saito, Lowest weights in cohomology of variations
of Hodge structure (preprint).

\item{[11]}
M.~Saito, Mixed Hodge modules, Publ. RIMS, Kyoto Univ.\ 26
(1990), 221--333.

\item{[12]}
M.~Saito, Admissible normal functions, J. Alg.\ Geom.\ 5 (1996),
235--276.

\item{[13]}
W.~Schmid, Variation of Hodge structure: The singularities of the
period mapping, Inv. Math. 22 (1973), 211--319.

\item{[14]}
J.H.M.~Steenbrink and S.~Zucker, Variation of mixed Hodge structure,
I, Inv.\ Math. 80 (1985), 489--542.

\item{[15]}
R.~Thomas, Nodes and the Hodge conjecture, J. Alg.\ Geom.\ 14
(2005), 177--185.

{\sfont
\medskip
RIMS Kyoto University, Kyoto 606-8502 Japan

\smallskip\vers}}
\bye